\documentclass[a4paper]{article}

\usepackage{amsfonts}
\usepackage{amsmath}
\usepackage{amsthm}\usepackage{amssymb}
\usepackage[mathscr]{eucal}
\usepackage{graphics}
\usepackage[cp866]{inputenc}
\usepackage[english]{babel}
\usepackage{makeidx}
\usepackage{latexsym,amsfonts,amssymb,amsmath,longtable,amsthm}
\input amssym.def
\usepackage{latexsym}
\usepackage{graphicx}
\usepackage[all]{xy}
\usepackage{amsfonts}
\usepackage{amsmath}
\usepackage{mathrsfs}
\usepackage{color}
\usepackage{ulem}
\usepackage[table]{xcolor}

\newtheorem{ttt}{Theorem}

\newtheorem{lem}{Lemma}
\newtheorem{cor}{Corollary}
\newtheorem{rem}{Remark}

\def\XXint#1#2#3{{\setbox0=\hbox{$#1{#2#3}{\int}$}
     \vcenter{\hbox{$#2#3$}}\kern-.5\wd0}}

 \usepackage[pagebackref,
            colorlinks,linkcolor={blue},citecolor={red}]
           {hyperref}

\def\const{{\mathrm{const}}}
\newcommand{\SA}{{\mathrm{Sect}}_\alpha}
\newcommand{\SAA}{{\mathrm{Sect}}_{\frac{\alpha}3}}

\def\0{{\mathbf{0}}}

\def\R{{\mathbb R}}
\def\th{{\theta}}

\def\Ti{{\mathscr{T}}}

\def\H{{\mathcal{H}}}

\def\ue{{\mathbf{u}}}
\def\loc{{\mathrm{loc}}}

\def\div{\hbox{\rm div}\,}
\def\arg{{\mathrm{arg}\,}}

\def\dist{\hbox{\rm dist}\,}

 \def\into{\int\limits_}

\def\u{\mathbf u}
\def\brf{\bar f}
\def\bu{\bar{\mathbf u}}

\def\e{\varepsilon}

\def\n{\mathbf n}

\def\div{\hbox{\rm div}\,}

\newcommand{\meas}{\mathop{\mathrm{meas}}}

\begin{document}
\title{On convergence of arbitrary $D$-solution of steady  Navier--Stokes system in $2D$ exterior domains\footnote{
Navier Stokes equations, two--dimensional exterior domains,
asymptotic behavior.}}

\author{ Mikhail V. Korobkov\footnote{School of Mathematical Sciences,
Fudan University, Shanghai 200433, China; and Novosibirsk State University,
1 Pirogova str., Novosibirsk, 630090, Russia;
korob@math.nsc.ru},  Konstantin Pileckas\footnote{Faculty of
Mathematics and Informatics, Vilnius University, Naugarduko Str.,
24, Vilnius, 03225  Lithuania; konstantinas.pileckas@mif.vu.lt} \,
and Remigio Russo\footnote{Dipartimento di Matematica e Fisica
Universit\`a degli studi della Campania "Luigi Vanvitelli," viale
Lincoln 5, 81100, Caserta, Italy; e-mail:
remigio.russo@unicampania.it}}

 \maketitle

 \date{}

\begin{abstract}
{We study solutions to
stationary Navier--Stokes system in two dimensional exterior
domain. We prove that any such solution with finite
Dirichlet integral converges to a constant vector at infinity uniformly. No additional condition (on symmetry or smallness, etc.) are assumed. The proofs based on arguments of the classical Amick's article  (Acta Math. 1988) and on results of a~recent paper by authors (arXiv 1711:02400) where the uniform boundedness of these solutions was established. }

\end{abstract}

\setcounter{section}{0}

\section{Introduction}

\setcounter{equation}{0}
\setcounter{ttt}{0}
\setcounter{lem}{0}

Let $\Omega$ be an exterior domain in $\R^2$, in particular,
\begin{equation}\label{Omega}
\Omega\supset{\mathbb R}^2\setminus B,
\end{equation}
where $B=B_{R_0}$ is the disk of radius $R_0$ centered at the origin with $\partial \Omega\subset B$.

We consider the solutions to the steady Navier--Stokes system
\begin{equation}
\label{SNS}
\left\{\begin{array}{r@{}l}
\nu \Delta{\u}-(\u\cdot\nabla)\u-\nabla p  & {} ={\bf 0}\qquad \hbox{\rm in } \Omega, \\[2pt]
\div{\u} & {} =0\,\qquad \hbox{\rm in } \Omega.
 \end{array}\right.
\end{equation} Starting from the pioneering
papers by J. Leray \cite{Ler}  it is now customary to consider
solutions to  (\ref{SNS}) with   finite Dirichlet integral
 \begin{equation}
\label{SxxS} \into\Omega|\nabla\u|^2<+\infty,
\end{equation}
known also as {\it D--solutions\/}. As is well known (e.g., \cite{OAL}), such solutions are
real--analytic in $\Omega$. The existence of solutions to \eqref{SNS} was also studied in \cite{FS}, \cite{KRMann}, \cite{kprplane}, \cite{AR}.

The problem  of the asymptotic behavior at infinity of an arbitrary
$D$--solution $(\u,p)$ to (\ref{SNS})  was
tackled by D.~Gilbarg \& H. Weinberger \cite{GW1}--\cite{GW2} and
Ch.~Amick \cite{Amick}. In \cite{GW2} it is shown that
\begin{equation}\label{prc}
p(z)-p_0=o(1)\qquad\mbox{ as }r\to\infty,
\end{equation}
i.e., pressure has a limit at infinity (one can choose, say, $p_0=0$\,) and
\begin{equation}
\label{SzzS}
\begin{array}{ l}
 \u(z) =o(\log^{1/2} r),\\[2pt]
 \omega(z)  =o(r^{-3/4}\log^{1/8}r),\\[2pt]
\nabla\u(z)=o(r^{-3/4}\log^{9/8}r),
 \end{array}
\end{equation}
where $r=|z|$ and
$$
\omega=\partial_2u_1-\partial_1u_2
$$
is the vorticity. If, in addition,  $\u$ is
bounded, then  there is a constant vector $\u_\infty$ such that
 \begin{equation}
\label{GGWW}
\displaystyle\lim_{r\to+\infty}\into0^{2\pi}|\u(r,\theta)-\u_\infty|^2d\theta=0,\end{equation}
and
\begin{equation}
\label{Sz123zS}
\begin{array}{l}
 \omega(z)  =o(r^{-3/4}),\\[2pt]
\nabla\u(z)=o(r^{-3/4}\log r).
 \end{array}
\end{equation}
Here if $\u_\infty={\bf 0}$, then
\begin{equation}
\label{as0}
\u(z)\to0\qquad\mbox{ uniformly as \ }|z|\to\infty.
\end{equation}
In the~case $\u_\infty\ne0$
  D.~Gilbarg \& H. Weinberger proved that there exists a~sequence of radii  $R_n\in(2^n,2^{n+1})$, $n\ge n_0$, such that
\begin{equation}
\label{GGWW-u9}
\displaystyle\sup\limits_{\theta\in[0,2\pi]}|\u(R_n,\theta)-\u_\infty|\to0\qquad\mbox{ as \ }n\to\infty.\end{equation}

In the classical and very elegant paper
\cite{Amick} Ch.Amick proved  that under zero boundary condition
\begin{equation}
\label{as-v}
\u|_{\partial\Omega}\equiv0
\end{equation}
the solution has the following asymptotic  properties:

\begin{itemize}
\item[(i)] \, $\u$ is bounded and, as a consequence, it satisfies
(\ref{GGWW}), (\ref{Sz123zS});

\item[(ii)] \,the
total head pressure $\Phi=p+\frac12|\u|^2$ and the absolute value
of the velocity $|\u|$ have the uniform limit at infinity, i.e.,
\begin{equation}
\label{dom6} |\u(r,\th)|\to|\u_\infty|\qquad\mbox{ as }\
r\to\infty,
\end{equation}
where $\u_\infty$ is the~constant vector from the
condition~(\ref{GGWW}).
\end{itemize}

Recently M.Korobkov, K.Pileckas and R.Russo~\cite{kpr-arxiv}  simplified the issue and proved that the first claim (i) holds in the~general case of $D$-solutions without (\ref{as-v}) assumption:

\begin{ttt}[\cite{kpr-arxiv}]
\label{T1} {\sl Let $\u$ be a $D$-solution to the Navier--Stokes
system~(\ref{SNS})
in the exterior domain~$\Omega\subset\R^2$.
Then $\u$ is uniformly bounded in $\Omega_0=\R^2\setminus B$, i.e.,
\begin{equation}
\label{dom2} \sup\limits_{z\in\Omega_0}|\u(z)|<\infty,
\end{equation}
where $B=B_{R_0}$ is an open disk with
sufficiently large radius: $B\supset\partial\Omega$.}
 \end{ttt}

Using the above--mentioned results of  D.~Gilbarg and
H.~Weinberger, we obtain immediately

\begin{cor}\label{GW-cor}
{\sl Let $\u$ be a $D$-solution to the Navier--Stokes
system~(\ref{SNS}) in a neighbourhood of infinity.  Then the
asymptotic properties~(\ref{prc}), (\ref{GGWW})--(\ref{Sz123zS})
hold.}
\end{cor}

 The  main result of the present paper is as follows.

\begin{ttt}
\label{T2} {\sl Let $\u$ be a $D$-solution to the Navier--Stokes
system~(\ref{SNS})
in the exterior domain~$\Omega\subset\R^2$.
Then $\u$ converges uniformly at infinity, i.e.,
\begin{equation}
\label{as-0}
\u(z)\to \u_\infty\qquad\mbox{ uniformly as \ }|z|\to\infty,
\end{equation}
where $\u_\infty\in\R^2$ is the~constant vector from the~equality~(\ref{GGWW}).}
 \end{ttt}

 The proof of Theorem~\ref{T2}  is based on a~combination of ideas of papers~
\cite{Amick}, \cite{kpr-arxiv} and \cite{GW2}.

If $\u_\infty\neq 0$, then by results of   L.I. Sazonov~\cite{Sazonov}, the convergence~(\ref{as-0}) ensures  that the solution behaves at
infinity as that of the linear Oseen equations (see also~\cite{Galdibook}).

\section{Notations and preliminaries}
\setcounter{equation}{0}

By {\it a domain} we mean an open connected set. We use standard
notations for Sobolev  spaces \,$W^{k,q}(\Omega)$, where $k\in{\mathbb N}$,
$q\in[1,+\infty]$. In our notation we do not
distinguish function spaces for scalar and vector valued
functions; it is clear from the context whether we use scalar or
vector (or tensor) valued function spaces.

For $q\ge1$ denote by $D^{k,q}(\Omega)$ the set of functions $f\in
W^{k,q}_{\loc}(\Omega)$ such that
$\|f\|_{D^{k,q}(\Omega)}=\|\nabla^k f\|_{L^q(\Omega)}<\infty$.

We denote by $\H^k$ the
$k$-dimensional Hausdorff measure, i.e.,
$\H^k(F)=\lim\limits_{t\to 0+}\H^k_t(F)$,
where $$\H^1_t(F)=\bigl(\frac{\alpha_k}2\bigr)^k\inf\{\sum\limits_{i=1}^\infty \bigl({\rm
diam} F_i\bigr)^k:\, {\rm diam} F_i\leq t, F\subset
\bigcup\limits_{i=1}^\infty F_i\}$$
and $\alpha_k$ is a Lebesgue volume of the unit ball in~$\R^k$.

In particular, for a curve $S$ the value $\H^1$ coincides with its length, and for sets $E\subset\R^2$ the $\H^2(E)$ coincides with the usual Lebesgue measure in~$\R^2$.

Also, for a curve $S$ by $\int\limits_Sf\,ds$ we denote the usual integral with respect to $1$-dimensional Hausdorff measure (=length). Further, for a set $E\subset\R^2$ by
$\int\limits_Ef(x)\,d\H^2$ or simply $\int\limits_Ef(x)$ we denote we integral with respect to the two-dimensional Lebesgue measure.

Below we present some  usual  results concerning the behaviour of~$D$-functions.

 \begin{lem}
\label{lt1} {\sl Let $f\in D^{1,2}(\Omega)$ and assume that
$$\int\limits_{D}|\nabla f|^2\,d\H^2<\e^2$$ for some $\e>0$ and for some ring
$D=\{z\in\R^2:r_1<|z-z_0|<r_2\,\}\subset\Omega$. Then the~estimate \begin{equation}
\label{est-1}
|\brf(r_2)-\brf(r_1)|\le \e\sqrt{\ln\frac{r_2}{r_1}}
\end{equation} holds, where $\brf$ means the mean value of~$f$ over the circle $S(z_0,r)$:
$$\brf(r):=\frac1{2\pi r}\int\limits_{|z-z_0|=r}f(z)\,ds.$$ }
\end{lem}
\textsc{Proof}.  Let $(r, \theta)$ be polar coordinates with the center in the point $z_0$. We have
$$
|\bar f(r_2)-\bar f(r_1)| = \Big|\int\limits_{r_1}^{r_2}\bar f'(r)dr\Big|\leq \int\limits_{r_1}^{r_2}\int\limits_0^{2\pi}\Big|\frac{\partial}{\partial r}f(r, \theta)\Big| d\theta dr\leq \int\limits_{r_1}^{r_2}\int\limits_0^{2\pi}\big|\nabla f(z)\big| d\theta dr.
$$
Estimating the  right-hand side by the Cauchy--Schwarz inequality we obtain
$$
|\bar f(r_2)-\bar f(r_1)|\leq \sqrt{\ln \frac{r_2}{r_1}}\Big( \int\limits_{r_1}^{r_2}\big(\int\limits_{|z-z_0|=r}\big|\nabla f(z)\big|^2ds\big) dr\Big)^{1/2}
$$
$$
\leq \sqrt{\ln \frac{r_2}{r_1}} \Big(\int\limits_{D}|\nabla f|^2\,d\H^2\Big)^{1/2}\leq  \varepsilon \sqrt{\ln \frac{r_2}{r_1}}.
$$
\hfill$\qed$

 \begin{lem}
\label{lt2} {\sl Fix a number $\beta\in(0,1)$. Let $f\in D^{1,2}(\Omega)$ and assume that
$$\int\limits_{D}|\nabla f|^2\,d\H^2<\e^2$$ for some $\e>0$ and for  some ring
$D=\{z\in\R^2:\beta R<|z-z_0|<R\,\}\subset\Omega$. Then there exists  a number $r\in[\beta R, R]$ such that
 the
estimate \begin{equation}
\label{est-2}
\sup\limits_{|z-z_0|=r}|f(z)-\brf(r)|\le c_\beta\e
\end{equation} holds, where the constant~$c_\beta$ depends on $\beta$ only.   }
 \end{lem}

\textsc{Proof} (see the proof of Lemma 2.2 in \cite{GW2}).  Take the polar coordinate system with the center at the~point~$z_0$. Since $\int\limits_{\beta R}^R\frac{1}{\rho}\int\limits_0^{2\pi} \big|\frac{\partial}{\partial \theta}f(\rho,\theta)\big|^2d\theta d\rho \leq \int\limits_D|\nabla f(z)|^2dz$, by the integral mean value theorem, there exists some $r\in [\beta R,R]$ such that
$$
 \int\limits_0^{2\pi} \bigl|\frac{\partial }{\partial\theta'}f(r,\theta')\bigr|^2d\theta  \leq \tilde c_\beta\int\limits_D|\nabla f(z)|^2dz.
$$
Therefore, by Holder inequality
\begin{equation}
\label{as-0-y}
 \int\limits_0^{2\pi} \big|\frac{\partial}{\partial \theta}f(r,\theta)\big|d\theta  \leq \biggl(2\pi  \int\limits_0^{2\pi}\bigl|\frac{\partial }{\partial\theta'}f(r,\theta')\bigr|^2d\theta\biggr)^{\frac12}\le c_\beta\e
\end{equation}
On the other hand,
 $$
 f(r,\theta)-f(r, \varphi) = \int\limits_\varphi^\theta\frac{\partial }{\partial\theta'}f(r,\theta')d\theta'.
 $$
Integrating this equality with respect to $\varphi$ and taking the average, we find
$$
|f(r,\theta)-\brf(r)|\le\int\limits_0^{2\pi}\bigl|\frac{\partial }{\partial\theta'}f(r,\theta')\bigr|\,d\theta'\le c_\beta\e.
$$
\hfill$\qed$

Summarize the results of these lemmas, we receive
\begin{lem}\label{lt3}
{\sl Under conditions of Lemma~\ref{lt2}, there exists $r\in [\beta R,R]$ such that
\begin{equation}
\label{est-3}
\sup\limits_{|z-z_0|=r}|f(z)-\brf(R)|\le \tilde c_\beta\e.
\end{equation}}

\end{lem}

\section{Proof of the main Theorem~\ref{T2}.}
\setcounter{equation}{0}
\setcounter{lem}{0}

Suppose the assumptions of Theorem~\ref{T2} are fulfilled. By
classical regularity results for $D$-solutions to the Navier--Stokes
system (e.g., \cite{Galdibook}),  the functions~$\u$ and $p$ are real--analytical on the set~$\Omega_0=\R^2\setminus B_{R_0}$.
Moreover, it follows from  results in~\cite{GW2} and Theorem~\ref{T1}, that  $\u$ and $p$ are  uniformly bounded in~$\Omega_0$,
\begin{equation}
\label{dom10} \sup\limits_{z\in\Omega_0}\bigl(|p(z)|+|\u(z)|\bigr)\le C<+\infty,
\end{equation}and the pressure $p$ has a limit at infinity; we could assume without loss of generality that
\begin{equation}
\label{asp-1}
p(z)\to 0\qquad\mbox{ uniformly as \ }|z|\to\infty.
\end{equation}
It is also well known (see \cite{Galdibook}) that all  derivatives of $\u$ uniformly converge to zero:
\begin{equation}
\label{asp-k}
\forall k=1,2,\dots\qquad\nabla^k\u(z)\to0\qquad\mbox{ uniformly as \ }|z|\to\infty.
\end{equation}
Further,  it is proved in~\cite{GW2} that  there exists a vector $\u_\infty\in \R^2$ such that
\begin{equation}
\label{GGWW-u1}
\displaystyle\lim_{r\to+\infty}\into0^{2\pi}|\u(r,\theta)-\u_\infty|^2d\theta=0,
\end{equation}
moreover, if $\u_\infty={\bf 0}$, then
\begin{equation}
\label{asp-0}
\u(z)\to0\qquad\mbox{ uniformly as \ }|z|\to\infty.
\end{equation}
Thus if $\u_\infty=0$, the statement of Theorem~\ref{T2} is known  and
we need to consider only the case
\begin{equation}
\label{asp-4}
\u_\infty\ne0.
\end{equation}

Consider the vorticity
$\omega=\partial_2u_1-\partial_1u_2$
which will play  the key role in our proof. Recall that~$\omega$  satisfies the elliptic equation
\begin{equation}
\label{eqv1}
\nu \Delta\omega=(\u\cdot\nabla)\omega.
\end{equation}
In particular, $\omega$ satisfies two-sided maximum principle in $\R^2$; moreover,
\begin{equation}
\label{eqv2}
\int\limits_{\Omega_0}r|\nabla\omega|^2<\infty
\end{equation}
(see~\cite{GW2}\,).

We will need also the following statement.

\begin{lem}
\label{LemmaGW}
{\sl Let $\u$ be a $D$-solution to the Navier--Stokes
system~(\ref{SNS})
in the exterior domain~$\Omega\subset\R^2$. Denoted by $\bu(z,r)$ the mean value  of $\u$ over the circle $S(z,r)$:
\begin{equation}
\label{eqv-c3}
\bu(z,r)=\frac1{2\pi r}\int\limits_{|\xi-z|=r}\u(\xi)\,ds
\end{equation}
and let $\varphi(z,r)$ be the argument of the complex number associated to the vector~$\bu(z,r)=(\bar u_1(r), \bar u_2(r))$, i.e., $\varphi(z,r)=\arg(\bar u_1(r)+i \bar u_2(r))$.
Suppose $|z|$ is large enough so that the disk $D_{z}=\bigl\{\xi\in \R^2:|\xi-z|\le\frac45|z|\bigr\}$ is contained in~$\Omega$. Assume also  that
$$|\bu(z,r)|\ge\sigma.$$
for some positive constant~$\sigma>0$ and for all $r\in\bigl(0,\frac45|z|\bigr]$. Then  the estimate
\begin{equation}
\label{eqv-c4--}
\sup\limits_{0< \rho_1\le\rho_2\le\frac45|z|}|\varphi(z,\rho_2)-\varphi(z,\rho_1)|\le \frac1{4\pi\sigma^2}\int\limits_{D_z}\biggl(\frac1r|\nabla\omega|+
|\nabla\u|^2\biggr)\,d\H_\xi^2
\end{equation}
holds, where $r=|\xi-z|$.
}
\end{lem}
For the proof of the estimate \eqref{eqv-c4--} see \cite[Proof of Theorem~4, page~399]{GW2}.\\

To apply the last Lemma~\ref{LemmaGW}, we need also the following simple technical assertion.

\begin{lem}
\label{LemmaGW-cor}
{\sl Let $\u$ be a $D$-solution to the Navier--Stokes
system~(\ref{SNS})
in the exterior domain~$\Omega\subset\R^2$. For $z\in\Omega$ denote as above
$$D_{z}=\bigl\{\xi\in \R^2:|\xi-z|\le\frac45|z|\bigr\}.$$
Then the uniform convergence
\begin{equation}
\label{eqv-c4}
\int\limits_{D_z}\frac1r|\nabla\omega|\,d\H_\xi^2\to 0\qquad\mbox{ as }|z|\to\infty
\end{equation}
holds, where again $r=|\xi-z|$.
}
\end{lem}

\textsc{Proof.} Take and fix arbitrary $\e>0$. Take also numbers $r_2>r_1>0$ large enough so that
\begin{equation}
\label{eqv-c9}
2\pi<\e r_1;
\end{equation}
\begin{equation}
\label{eqv-c10}
\int\limits_{D_z}r|\nabla\omega|^2\,d\H_\xi^2<\e\qquad\mbox{ if }|z|>r_2;
\end{equation}
\medskip
\begin{equation}
\label{eqv-c11}
\qquad\qquad 2\pi r_1\max\limits_{|\xi-z|<r_1}|\nabla\omega(\xi)|<\e\qquad\mbox{ if }|z|>r_2
\end{equation}
(the existence of such numbers follows from the estimate~(\ref{eqv2}) and from the uniform convergence~(\ref{asp-k})\,).

Now take  arbitrary $z\in\R^2$ with $|z|>r_2$.
Then the disk $D_z$ is represented as the union of two sets $D_z=D_1\cup D_2$, where
$$D_1=\bigl\{\xi\in\R^2:|\xi-z|<r_1\bigr\},\qquad D_2=\bigl\{\xi\in\R^2: r_1\le|\xi-z|<\frac45|z|\bigr\}.$$
We have
$$
\int\limits_{D_1}\frac1r|\nabla\omega|\,d\H_\xi^2<\max\limits_{|\xi-z|<r_1}|\nabla\omega(\xi)|\int\limits_{D_1}\frac1r\,d\H_\xi^2
$$
\begin{equation}
\label{eqv-c12}
=
 2\pi r_1\max\limits_{|\xi-z|<r_1}|\nabla\omega(\xi)|\overset{\mbox{\footnotesize (\ref{eqv-c11})}}<\e.
\end{equation}
Further, applying the elementary inequality
$\frac1r|\nabla\omega|<\frac1{r^3}+r|\nabla\omega|^2$, for the domain $D_2$ we have:
$$
\int\limits_{D_2}\frac1r|\nabla\omega|\,d\H_\xi^2<\int\limits_{D_2}\frac1{r^3}\,d\H_\xi^2+\int\limits_{D_2}r|\nabla\omega|^2\,d\H_\xi^2
$$
\begin{equation}
\label{eqv-c13}
=
2\pi\int\limits_{r=r_1}^{\frac45|z|}\frac1{r^2}\,dr+\int\limits_{D_2}r|\nabla\omega|^2\,d\H_\xi^2\overset{\mbox{\footnotesize(\ref{eqv-c9})-(\ref{eqv-c10})}}<2\e.
\end{equation}
From the inequalities~(\ref{eqv-c12})--(\ref{eqv-c13})  it follows that
\begin{equation}
\label{eqv-c15}
\int\limits_{D_z}\frac1r|\nabla\omega|\,d\H_\xi^2<3\e.
\end{equation}
We proved the last inequality for any~$z\in\R^2$ with $|z|>r_2$. Since the number $\e>0$ is arbitrary, the required
convergence~(\ref{eqv-c4}) is established. \hfill $\qed$

\

Further we will use the following two criteria for the uniform convergence of the velocity:

\begin{lem}
\label{T3} {\sl Let $\u$ be a $D$-solution to the Navier--Stokes
system~(\ref{SNS})
in the exterior domain~$\Omega\subset\R^2$.
Suppose that at least one of the following two conditions is fulfilled:

\begin{itemize}
\item[{\rm(i)}] \,$\omega(z)=o(|z|^{-1})$ \ as $|z|\to\infty$;

\item[{\rm(ii)}] \,the absolute value of the velocity has a uniform limit at infinity:
\begin{equation}
\label{asp-11}
|\u(z)|\to|\u_\infty|\qquad\mbox{ uniformly as \ }|z|\to\infty,
\end{equation}
where the vector $\u_\infty$ was specified above.\end{itemize}
Then $\u$ converges uniformly at infinity as well, i.e., the formula~(\ref{as-0}) holds. }
 \end{lem}

\textsc{Proof.}
 Part~(i) was established by Amick  (see \cite{Amick}, Remark~3(i) on p.~103 and the proof of Theorem~19).
Recall, that his argument is based on the classical  Cauchy-type representation formula of complex analysis:
\begin{equation}
\label{eqv-c1}
w(z)=\frac1{2\pi i}\oint\limits_{|\xi-z_0|=r}\frac{w(\xi)\,d\xi}{\xi-z}+\frac1{2\pi i}\iint\limits_{|\xi-z_{0}|<r}\frac{\omega(\xi)}{\xi-z_0}\,dx\,dy,
\end{equation}
where $w(\xi)=u_1(\xi)-iu_2(\xi)$ \,and  \,$\xi=x+iy$.

\medskip Let us prove the second part of Lemma~\ref{T3}. Suppose that assumption~(ii) is fulfilled. If $\u_\infty=0$, then there is nothing to prove (see the above discussion concerning the results of  D.~Gilbarg \& H. Weinberger \cite{GW1}--\cite{GW2}\,). So we assume without loss of generality that
\begin{equation}
\label{exp0-}
|\u_\infty|>0.
\end{equation}
From assumption~(\ref{asp-11}) and  Lemmas~\ref{lt1}--\ref{lt3} it follows that
\begin{equation}
\label{exp0}
\sup\limits_{0< \rho\le\frac45|z|}\biggl|\,|\u_\infty|-|\bu(z,\rho)|\biggr|\to 0\qquad\mbox{ uniformly as \ }|z|\to\infty,
\end{equation}
where  $\bu(z,r)$ is the mean value of $\u$ over the circle $S(z,r)$.
In particular, because of inequality~(\ref{exp0-}),
 there exist numbers $\sigma>0$ and $R_*>0$ such that
\begin{equation}
\label{eqv-c2}
|\bu(z,r)|\ge \sigma\qquad\mbox{ if }\ |z|\ge R_*\quad\mbox{ and }\ 0< r\le\frac45 |z|.
\end{equation}
 Then, by Lemma \ref{LemmaGW}, the argument $\varphi (z,r) $ of the complex number associated to $\bu(z,r)$ satisfies the estimate \eqref{eqv-c4--}.
From  \eqref{eqv-c4--}--\eqref{eqv-c4} \ it follows immediately that
\begin{equation}
\label{eqv-c5}
\sup\limits_{0< \rho_1\le\rho_2\le\frac45|z|}|\varphi(z,\rho_2)-\varphi(z,\rho_1)|\to 0
\end{equation}
uniformly as~$|z|\to\infty$.
In particular,
\begin{equation}
\label{exp1}
\sup\limits_{0< \rho\le\frac45|z|}|\arg\u(z)-\arg\bu(z,\rho)|\to 0
\end{equation}
uniformly as~$|z|\to\infty$.
From the assumption  (\ref{asp-11}) and (\ref{exp0}) we have
\begin{equation}
\label{exp00---}
\sup\limits_{0< \rho\le\frac45|z|}\biggl|\,|\u(z)|-|\bu(z,\rho)|\biggr|\to 0\qquad\mbox{ uniformly as \ }|z|\to\infty.
\end{equation}
Summarizing the information from formulas~(\ref{exp1})--(\ref{exp00---}),
we obtain
\begin{equation}
\label{exp3}
\sup\limits_{0< \rho\le\frac45|z|}|\u(z)-\bu(z,\rho)|\to 0\qquad\mbox{ uniformly as \ }|z|\to\infty.
\end{equation}

Consider the sequence of  circles
$S_{R_n}=\{\xi\in \R^2:|\xi|=R_n\}$ such that $2^n<R_n<2^{n+1}$ and
\begin{equation}
\label{GGWW-u9--}
\displaystyle\sup\limits_{|\xi|=R_n}|\u(\xi)-\u_\infty|=\e_n\to0\qquad\mbox{ as \ }n\to\infty
\end{equation}
(the existence of such sequence is guaranteed by above mentioned results of  D.~Gilbarg and H. Weinberger, see ~(\ref{GGWW-u9})\,).

Now take a point $z\in\R^2$ with sufficiently large~$|z|$ and take also the natural number~$n=n_z$ such that
$$2^{n+1}\le |z|<2^{n+2}.$$
Then by construction and by the triangle inequality we have
\begin{equation}
\label{exp--1}
S_{R_n}\cap S_{z,\rho}\ne\emptyset \qquad\mbox{ if }\quad\frac34|z|<\rho<\frac45|z|,
\end{equation}
where $S_{z,\rho}=\{\xi\in\R^2:|\xi-z|=\rho\}.$  From Lemma~\ref{lt2} it follows that there exists $\rho_*\in\bigl(\,\frac34|z|,\frac45|z|\,\bigr)$ such that
\begin{equation}
\label{exp7}
\sup\limits_{|\xi-z|=\rho_*}|\u(\xi)-\bu(z,\rho_*)|=\e_z,
\end{equation}
where $\e_z\to0$ \  uniformly as \ $|z|\to\infty$. Summarizing the information from formulas~(\ref{GGWW-u9--})--(\ref{exp7}), we obtain that
\begin{equation}
\label{exp8}
|\u_\infty-\bu(z,\rho_*)|=\e'_z\to0\qquad\mbox{ uniformly as \ }|z|\to\infty.
\end{equation}
Finally, from the last formula and from (\ref{exp3}) we conclude that
\begin{equation}
\label{exp00}
|\u_\infty-\u(z)|\to0\qquad\mbox{ uniformly as \ }|z|\to\infty,
\end{equation}
as required.
 The Lemma~\ref{T3} is proved completely. \hfill $\qed$\\

\textsc{Proof of Theorem~\ref{T2}}.
For a point $z\in\Omega_0$ denote by $K(z)$ the connected component of the level set of the vorticity~$\omega$ containing~$z$, i.e.,
$K(z)\subset\{x\in\Omega_0:\omega(x)=\omega(z)\}$. Here we understand the notion of connectedness in the sense of general topology.

We consider two possible cases:

{\bf Case I}.  \textsf{Level sets of $\omega$ separate infinity from the origin:}
\begin{equation}
\label{case1}
\exists z_*\in \Omega_0:\ \ \omega(z_*)\ne 0\ \ \mbox{ and }\ \ K(z_*)\cap\partial\Omega_0=\emptyset.
\end{equation}

{\bf Case II}.   \textsf{Level sets of $\omega$ do not separate infinity from the origin:}
\begin{equation}
\label{case2}
K(z)\cap\partial\Omega_0\ne\emptyset\qquad \forall z\in \Omega_0,
\end{equation}

In Case I, we shall show that
\begin{equation}\label{case11}
|z|\omega(z)\to 0  \qquad\mbox{ uniformly as \ }|z|\to\infty
\end{equation}
and we obtain the statement of Theorem applying Lemma \ref{T3}(i).

In Case II, we prove that
\begin{equation}\label{case21}
|\u(z)|\to|\u_\infty|\qquad\mbox{ uniformly as \ }|z|\to\infty,
\end{equation}
where  $\u_\infty$ is the vector defined in \eqref{GGWW-u1}. In this case the statement of Theorem will follow from Lemma \ref{T3}(ii). \\

Consider the case \eqref{case1}.  Note that then the set $K(z_*)$ is compact. Indeed, the set $K(z_*)$  is connected and if it is not compact, it should "reach" infinity. Since the vorticity  tends to zero at infinity, $\omega(z)$ has to be zero on $K(z_*)$, but this contradicts the assumption \eqref{case1}.

Next, by elementary compactness and continuity arguments we have that there exists $\delta_0>0$ such that
\begin{equation}\label{case21-exp1}
K(z)\mbox{ is a compact set satisfying \ $K(z)\cap\partial\Omega_0=\emptyset$ \ whenever \ $|z-z_*|<\delta_0$}.
\end{equation}
Note, that since $\omega$ is an analytical nonconstant function, we have that
$\omega(z)\ne\const$ in any open neighborhood of~$z_*$.

Recall, that a real  number $t$ is called {\it a regular value of~}$\omega$, if the set $\{z\in\Omega_0:\omega(z)=t\}$ is nonempty and
$\nabla\omega(z)\ne0$ whenever $\omega(z)=t$. By the classical Morse--Sard theorem, almost all values of~$\omega$ are regular.
Now take a point $z_1$ satisfying
$|z_1-z_*|<\delta_0$ with regular value~$t_1=\omega(z_1)$. Then by definition and regularity assumptions  the set $K(z_1)$ is a smooth compact curve (={\sl ``compact one dimensional
manifold without boundary''}).
By obvious topological reasons, $K(z_1)$ is a smooth curve homeomorphic to the circle. Since $\omega $ satisfies  maximum principle, this circle surrounds the origin. Therefore, the curve $K(z_1)$ separates the boundary~$\partial\Omega_0$ from infinity\footnote{It means that infinity and the set $\partial\Omega_0$ \ lie  in the different connected components of the~set $\R^2\setminus K(z_1)$.}.

Denote $R_*=\max\{|z|:z\in K(z_1)\}$ \ and \ $\Omega_*=\{z\in\R^2:|z|> R_*\}$. Then by construction
 we have
\begin{equation}
\label{1c-2}
 K(z)\cap\partial\Omega_0=\emptyset\qquad\forall z\in\Omega_*.
\end{equation}
Applying again the same  Morse--Sard theorem, we obtain  that for almost all $t\in\R\setminus\{0\}$
if $z\in\Omega_*$ and $\omega(z)=t$, then $K(z)$ is a smooth curve homeomorphic to the circle.
Since $\omega $ satisfies  maximum principle, we conclude  that this circle surrounds the origin, moreover,
\begin{equation}
\label{1c-3}
K(z_1)=K(z_2)\qquad\mbox{if } z_1,z_2\in\Omega_*\quad\mbox{ and }\omega(z_1)=\omega(z_2)\ne0.
\end{equation}
This implies that
\begin{equation}
\label{1c-4}
\omega(z) \mbox{\ \ does not change sign in }\Omega_*.
\end{equation}
Indeed, let there are points $z_1,z_2\in\Omega_*$ with regular values  $\omega(z_1)<0$ and $\omega(z_2)>0$. Taking into account that $\omega(z)$ is vanishing at the infinity, by maximum principle, $\omega(z)$ is negative in the exterior of $K(z_1)$ and $\omega(z)$ is positive in the exterior of $K(z_2)$. Since this is impossible, $\omega(z)$ cannot change the sign.

Thus we may suppose without loss of generality that
 \begin{equation}
\label{1c-5}
\omega(z)\ge 0\qquad\mbox{ in } \Omega_*.
\end{equation}
Then by the maximum principle we have the strict inequality
 \begin{equation}
\label{1c-5-s}
\omega(z)>0\qquad\mbox{ in }\Omega_*.
\end{equation}
Moreover, from (\ref{1c-3}) and  from the uniform convergence (see \eqref{asp-k})
 \begin{equation}
\label{1c-5--s}
\omega(z)\to0\qquad\mbox{ as }|z|\to\infty
\end{equation}
and from  Morse--Sard theorem   we conclude that
there exists a number $\delta>0$  such that
\begin{equation}
\label{1c-6}
\begin{array}{ l}
\mbox{for almost all $t\in (0,\delta)$ the set $K_t:=\{z\in\Omega_*:\omega(z)=t\}$}\\[2pt]
\mbox{coincides with the smooth curve homeomorphic to the circle}\\[2pt]
\mbox{such that $K_t\cap\partial\Omega_*=\emptyset$ and $\nabla\omega\ne0$ on $K_t$}.
 \end{array}
\end{equation}

Denote by $\Ti$  the set of full measure in the interval $(0,\delta)$ consisting of values~$t$ satisfying~(\ref{1c-6}).
Denote also by $\Omega_t$ the unbounded connected component of the set $\R^2\setminus K_t$.
Since $\omega$ satisfies the maximum principle,  the sets $K_t$ have the following monotonicity property:
\begin{equation}
\label{1c-7}
\Omega_{t_1}\subset\Omega_{t_2}\qquad\mbox{ if }\ \ 0<t_1<t_2.
\end{equation} Moreover, from the uniform convergence (\ref{1c-5--s}), it follows that
\begin{equation}
\label{1c-8}
\inf\{|z|:z\in\Omega_t\}\to \infty\qquad\mbox{ as }\ \ t\to0+.
\end{equation}

Our task is to show the property (i) of Lemma~\ref{T3}, i.e., to show that
\begin{equation}
\label{1c-9}
|z|\omega(z)\to0 \qquad \mbox{uniformly as  }\  \ |z|\to\infty.
\end{equation}
The last condition is equivalent to
\begin{equation}
\label{1c10}
t g(t)\to0\qquad\mbox{ as }\ \ t\to0+,
\end{equation}
where the function $g(t)$ is defined by
\begin{equation}
\label{1c-11}
g(t):=\sup\{|z|:z\in K_t\}.
\end{equation}
Obviously, $g(t)\le \H^1(K_t)$, where, recall, $\H^1$ is the one-dimensional Hausdorff measure (=length).

For $t\in\Ti$ and $R>R_*$ denote $\Omega_{t,R}=\Omega_t\cap B_R=\{z\in\Omega_t:|z|<R\}$.  Then for sufficiently large $R$
$$
\partial\Omega_{t,R}=K_t \cup S_R,
$$
where $S_R=\{z\in\R^2:|z|=R\}$ is the corresponding circle.
Integrating the equation \eqref{eqv1}
over the domain~$\Omega_{t,R}$ \ and taking into account that $(\u\cdot\nabla)\omega=\div(\u \omega)$,  we obtain
 \begin{equation}
\label{1c-13}
\int\limits_{K_t}|\nabla\omega|\,ds+\int\limits_{S_R}\nabla\omega\cdot{\bf n}\,ds=t\int\limits_{K_t}\u\cdot\n\,ds+\int\limits_{S_R}\omega\,\u\cdot\n\,ds.
\end{equation}
Here ${\bf n}$ is a unit vector  of the outward with respect to $\Omega_{t,R}$ normal to  $\partial\Omega_{t,R}$. Note also that the unit normal to the level set $K_t=\{z\in \Omega_*: \omega(z)=t\}$ is given by the formula ${\bf n}=\dfrac{\nabla \omega}{|\nabla\omega|}$.

Since $\div\u = 0$, we have $\int\limits_{K_t}\u\cdot\n\,ds=\int\limits_{\partial\Omega_*}\u\cdot\n\,ds=C_*$, i.e., this value does not  depend on~$t$. On the other hand, the estimate $\int\limits_{\Omega_0}\bigl(|\omega|^2+|\nabla\omega|^2\bigr)\,d\H^2<\infty $ implies that there is a~sequence $R_k\to+\infty$ such that
$$\int\limits_{S_{R_k}}\bigl(|\omega|+|\nabla\omega|\bigr)\,ds\to0.$$
 Taking  $R=R_k$ in the  equality (\ref{1c-13}) and having in mind the uniform boundedness  of  the velocity (see~(\ref{dom2})\,), we deduce, passing $R_k\to+\infty$, that
 \begin{equation}
\label{1c-14}
\int\limits_{K_t}|\nabla\omega|\,ds=C_*t.
\end{equation}

Further, for $t\in(0,\frac12\delta)$ denote $E_t=\{z\in\Omega_*:\omega(z)\in(t,2t)\}.$ By construction,
$$\partial E_t=K_t\cup K_{2t}.$$
Applying the classical Coarea formula  (see, e.g.,~\cite{Maly})
$$
\int\limits_{E_t}f\,|\nabla \omega|\,d\H^2=\int\limits_t^{2t}\biggl(\int\limits_{K_\tau}f\,ds\biggr)\,d\tau
$$
 for $f=|\nabla\omega|$ we obtain
 \begin{equation}
\label{1c-15}
\int\limits_{E_t}|\nabla\omega|^2\,d\H^2=\int\limits_t^{2t}\biggl(\int\limits_{K_\tau}|\nabla\omega|\,ds\biggr)\,d\tau\overset{\footnotesize{(\ref{1c-14})}}=
\int\limits_t^{2t}C_*\tau\,d\tau=3 C_*t^2.
\end{equation}
Applying now the same Coarea formula for
$f=1$ and using the~Cauchy--Schwarz inequality, we get
\begin{equation}
\label{dom18}
\begin{array}{lcr}
\int\limits_t^{2t}\H^1(K_\tau)\,d\tau=\int\limits_{E_t}|\nabla
\omega|\,d\H^2\le\biggl(\int\limits_{E_t}|\nabla\omega|^2\,d\H^2\biggr)^{\frac12}
\biggl(\meas E_t\biggr)^{\frac12}\\
\\
\overset{\footnotesize{(\ref{1c-15})}}=\sqrt{3C_*}\biggr(t^2\meas(E_t)\biggr)^{\frac12}\le\sqrt{\frac34C_*}\biggr(\int\limits_{E_t}\omega^2\,d\H^2
\biggr)^{\frac12}\le\e_t\to0\quad\mbox{ as }t\to0.
\end{array}
\end{equation}
Here  we have used also the fact that $t\le|\omega(z)|\le2t $ in $E_t$.
By virtue of the  mean-value theorem, this implies that for any sufficiently small $t\in\Ti$ there exists a number $\tau\in[t,2t]$ such that
$$
t\H^1(K_\tau)\le\e_t.
$$
By construction, the closed curve $K_\tau$ surrounds~$K_{2t}$.  Therefore,
$$
\sup\{|z|:z\in K_{2t}\}\le \H^1(K_\tau)\le\frac{\e_t}{t}
$$
with $\e_t\to0$ as $t\to0$. From the last inequality we receive the relation~(\ref{1c10})
which is equivalent to~\eqref{1c-9}.  According to Lemma~\ref{T3}(i), this finishes the proof of  Theorem~\ref{T2} in the considered Case~I.

\

Consider Case II, i.e., the  when
\begin{equation}
\label{1case2}
K(z)\cap\partial\Omega_0\ne\emptyset\qquad \forall z\in \Omega_0.
\end{equation}
 Now we shall prove that the assertion \eqref{case21} is valid.

Let us recall that Ch.~Amick~\cite{Amick} has proved the convergence~(\ref{case21}) under the assumption that
\begin{equation}
\label{2c-asp1}
\u|_{\partial\Omega}=0.
\end{equation}
The condition~(\ref{case21}) was used in~\cite{Amick}  in order to define  the stream function~$\psi$ in the neighborhood
of infinity:
\begin{equation}
\label{exp-app}
\nabla\psi=\u^{\bot}=(-v,u),
\end{equation}
where $\u=(u,v)$. Using the stream function~$\psi$, Amick introduced an~auxiliary function~$\gamma=\Phi-\omega\psi$, where $\Phi:=p+\frac12|\u|^2$ is the~Bernoulli pressure.
The gradient of this  auxiliary function~$\gamma$ satisfies the identity
$$\nabla\gamma=-\nu \nabla^\bot\omega-\psi\nabla\omega.$$
Then $\nabla\gamma\cdot\nabla^\bot\omega=-\nu|\nabla^\bot\omega|^2$,
and therefore, $\gamma$ has the~following   monotonicity properties:
\begin{equation}
\label{exp-app2}
\begin{array}{lcr}
\mbox{$\gamma$ is monotone along
 level sets of the~vorticity~$\omega=c$ and }\\ [10pt]
\mbox{vice versa~--- the~vorticity $\omega$ is monotone along level sets of~$\gamma=c$,}
\end{array}
\end{equation}
see \cite{Amick}.

Obviously,
the stream function~$\psi$ (and, consequently, the corresponding
auxiliary function~$\gamma$\,) is well defined in the
neighborhood of infinity under the more general condition
\begin{equation}
\label{dom5} \int\limits_{\partial\Omega_0}\u\cdot \n\,ds=0
\end{equation}
instead of~(\ref{2c-asp1}).
However,  in the general case the flow-rate of the velocity field is not zero,
\begin{equation}
\label{dom5--} \int\limits_{\partial\Omega_0}\u\cdot \n\,ds\ne0,
\end{equation}
and, therefore, the stream  function $\psi$  can not be defined in the neighborhood of infinity.

 We will overcome this difficulty using the assumption~(\ref{1case2}).
Take and fix a radius $R_*>R_0$ \ ($R_*$ could be chosen arbitrary  large\,) and consider the domain~$\Omega_*=\{
z\in\R^2:|z|>R_*\}$. Denote by $U_i$ the connected components of the open set~$\{z\in\Omega_*:\omega(z)\ne0\}$. Then there holds the following

\begin{lem}\label{lll1}{\sl Under assumption~(\ref{1case2}) the following assertions are fulfilled:

\begin{itemize}

\item[(i)] \,There are only finitely many components~$U_i$, $i=1,\dots,N$;

\item[(ii)] \,Every $U_i$ is a simply connected open set;

\item[(iii)] \,The vorticity $\omega(z)$ change sign in every  neighborhood of infinity, i.e., there exist two sequences of points $z^+_n$ and $z^-_n$ such that
$\omega(z^+_n)>0$, \ $\omega(z^-_n)<0$  and $\lim\limits_{n\to\infty}|z^+_n|=\lim\limits_{n\to\infty}|z^-_n|=\infty$.
\end{itemize}}
\end{lem}

We shall prove Lemma~\ref{lll1}  below.  Let us finish the proof of the theorem using this lemma. The~components $U_i$ play
 also an  important role in the arguments of Amick.   In particular, he proves in~\cite{Amick} the same properties~(i)--(iii)
 using the boundary condition $\ue|_{\partial\Omega}= 0$.
 Here, in Lemma~\ref{lll1}, we get the properties~(i)--(iii) because of the assumption~(\ref{1case2}).
 Since $U_i$ are simply connected, this allows us to define the stream function~$\psi$ in every component~$U_i$. Moreover, since $\omega= 0$ on $\Omega_*\cap \partial  U_i$, the auxiliary function~$\gamma=\Phi-\omega\psi$ is well defined and continuous on the whole domain~$\Omega_*$.
 After the functions $\psi$ and $\gamma$ are defined, we can repeat the arguments of the paper~\cite{Amick} and to prove the convergence~(\ref{case21}) of absolute value of the velocity at infinity. By Lemma \ref{T3}(ii) this implies the statement  of Theorem~\ref{T2}. For the reader convenience
we recall  the corresponding arguments of Amick~\cite{Amick}  in Appendix (we also simplify some of his proofs). \hfill$\qed$

\textsc{Proof of Lemma \ref{lll1}.} Let us prove~(iii) first. Suppose this is not true, i.e.,
there exists $R_1>0$ such that $\omega(z)$ does not change sign in $\Omega_1=\{z:|z|>R_1\}$. Without loss of generality assume that $\omega(z)\ge0$ in~$\Omega _1$. Then, by maximum principle,
\begin{equation}\label{cas1} \omega(z)>0\qquad\mbox{ in }\Omega_1.
\end{equation}
Take arbitrary $R_2>R_1$ and denote
\begin{equation}\label{cas2} \delta:=\inf\limits_{z\in S_{R_2}}\omega(z),
\end{equation}
where, recall, $S_{R_2}=\{z\in\R^2:|z|=R_2\}$. \ By~(\ref{cas1}), $\delta>0$. Now take any  $z_2$ such that $|z_2|>R_2$ and $\omega(z_2)<\delta$. Then by construction
$K(z_2)\cap S_{R_2}=\emptyset$. Therefore, $K(z_2)\cap S_{R_0}=K(z_2)\cap\partial\Omega_0=\emptyset$, a contradiction with~(\ref{1case2}).

\medskip

 (ii).$\quad$ Fix a component~$U_i$ and take an arbitrary curve $S\subset U_i$ homeomorphic to the unit circle. By construction, there exists $\delta>0$ such that
$$\omega(z)>\delta\qquad\forall z\in S.$$
The curve $S$ split the plane $\R^2$ into the two components: $\R^2\setminus S=\Omega_S\cup \Omega_{\infty}$, where $\partial\Omega_S=\partial\Omega_\infty=S$, \
$\Omega_S$ is a bounded domain homeomorphic to the disk, and  $\Omega_{\infty}$ is a neighborhood of infinity. Now we have to consider two cases:

\begin{itemize}

\item[($\alpha$)] \,the curve $S$ surrounds the origin. Then $\Omega_\infty\subset\Omega_*$, and, by maximum principle, $\omega\ge 0$ in $\Omega_\infty$. Thus, we received the contradiction with property~(iii) proved just above.

\item[($\alpha\alpha$)] \,the curve $S$ does not surround the origin. Then $\Omega_S\subset\Omega_*$, and, by maximum principle, $\omega> 0$ in $\Omega_S$.  Therefore,
$\Omega_S\subset U_i$.  Since $S$ was arbitrary, it means that $U_i$ is a simply connected set.
\end{itemize}

\medskip Let us prove (i).  Since $\omega$ is a nonzero analytical function, the set $Z_*=\{z\in S_{R_*}:\omega(z)=0\}$ is finite (recall, that $S_{R_*}$ is a circle of radius~$R_*$\,). Let $S_j$, $j=1,\dots,M$, be
the connected components of the set $S_{R_*}\setminus Z_*$.

Fix arbitrary component~$U_i$.
By maximum principle, $\omega(z)$ is  not identically zero on~$\partial U_i$, i.e., there exists a point~$z_0$ such that $$z_0\in \partial U_i\qquad\mbox{ and }\qquad\omega(z_0)\ne0.$$
On the other hand, by definition~$U_i$ is a connected component of the open set
$$\{z\in\Omega_*:\omega(z)\ne0\},$$
in particular,
we have the identity $\omega(z)\equiv 0$ on  the set $\Omega_*\cap\partial U_i$. Therefore,
$$z_0\in \partial \Omega_*=S_{R_*}.$$
It means, using the above notation, that there exists a number $j(i)\in\{1,\dots,M\}$ such that
$$z_0\in S_{j(i)}.$$
Then by
 elementary properties of connected sets and by definitions of $S_j$ and $U_i$, we have
$$S_{j(i)}\subset\partial U_i,$$ and $$\biggl[\,j(i_1)=j(i_2)\,\biggr]\Rightarrow U_{ i_1}=U_{i_2},$$
i.e., the function $i\mapsto j(i)$ is injective. Finally, since the family of components $S_j$ is finite, we conclude that the family~$U_i$ is finite as well.
 This finishes the proof of Lemma~\ref{lll1}.  \hfill$\qed$

\

\section{Appendix}
\setcounter{equation}{0}
\setcounter{ttt}{0}
\setcounter{lem}{0}
\setcounter{cor}{0}

For reader's convenience  we  recall here some steps of the corresponding arguments of Amick~\cite{Amick} for the proof of the~convergence~(\ref{case21}).

Our Lemma~\ref{lll1} implies, in particular,  that
there exists at least one unbounded component $U_{k_1}$  where $\omega$ is strictly positive and at least one unbounded component $U_{k_2}$  where $\omega$ is strictly negative (cf. with~\cite[Theorem~8, page~84]{Amick}\,).

\medskip

First  of all we mention, that by
\cite[Theorem~15, page~95]{Amick}, if we take  the number ~$R_*$  large enough, then there holds the following statement
\begin{equation}\label{cas3} \nabla\omega(z) \ne0\quad\mbox{ if }\omega(z)=0\mbox{ \  and \ }|z|\ge R_*.
\end{equation}
This gives the possibility to clarify the geometrical and topological structure of the components~$U_i$.  Namely,
$\Omega_*\cap\partial U_i$ consists of finitely many smooth (even analytical) curves.

Let $U_i$, $i=1,\dots,M$ be a family of {\it unbounded} components~$U_i$.  Then Amick proved the following geometrical and analytical characterization for~them:

\begin{ttt}[see Theorem~11, page 89 in~\cite{Amick}]\label{ttt1}{\sl For every $U_i$, $i=1,\dots, M$,
\begin{itemize}

\item[($\alpha$)] \,The set $\Omega_*\cap \partial U_i$ has precisely two unbounded components which may be parametrised as $\{(x_j(s), y_j(s)):s\in (0,\infty)\}$, $j=1,2$. In addition, $(x_j(0), y_j(0))\in\{|z|=R_*\}$, $s$ denotes the arc-length measure from these points, and the functions $x_j(\cdot)$ and $y_j(\cdot)$ are real-analytical (if we choose $R_*$ large enough to have~{\rm(}\ref{cas3}{\rm)}\,). The function~$\omega$ vanishes on these arcs and $|(x_j(s), y_j(s))|\to\infty$ as $s\to\infty$.

\item[($\alpha\alpha$)] \,The maps $s\mapsto \Phi(x_j(s),y_j(s))$ are monotone decreasing and increasing on $(0,\infty)$, respectively, for $j=1$ and $j=2$.
\end{itemize}
}
\end{ttt}

Since the Bernoulli pressure $\Phi$ is uniformly bounded, by Weierstrass Monotone convergence theorem we have that the functions~$s\mapsto \Phi(x_j(s),y_j(s))$ have some limits as $s\to\infty$ for $j=1,2$. After the usual agreement that
\begin{equation}\label{cas4}p(z)\to 0\qquad\mbox{ as }|z|\to\infty,
\end{equation}
and taking into account the convergence on the family of circles~(\ref{GGWW-u9}) we obtain

\begin{cor}\label{ttt2}{\sl  Functions from item~($\alpha\alpha$) of Theorem~\ref{ttt1} have the same limit
\begin{equation}\label{cas5}\Phi(x_j(s),y_j(s))\to\frac12|\u_\infty|^2\qquad\mbox{ as }s\to\infty.
\end{equation}
}
\end{cor}

The next step concerns the auxiliary function~$\gamma$.
One of the most important tool in~\cite{Amick} is the following assertion.

\begin{ttt}[see Theorem~14, page 92 in~\cite{Amick}]\label{ttt3}{\sl For every $U_i$, $i=1,\dots, M$, the convergence
\begin{equation}\label{cas6}\gamma(z)\to \frac12|\u_\infty|^2\qquad\mbox{ uniformly as \ }|z|\to\infty, \ z\in U_i
\end{equation}
holds.
}
\end{ttt}

{\sc Proof}.  We reproduce  here a simplified version of the proof of Theorem~14 in  \cite[pages~92--94]{Amick}.

Take and fix an unbounded  component~$U_i$.  We  assume without loss of generality that $|\u_\infty|=1$\, and\, $\omega(z)>0$ in $U_i$. By construction, we have
\begin{equation}\label{expl-ap1}\mbox{$\omega\equiv0$ and $\gamma\equiv\Phi$ on $\Omega_*\cap\partial U_i$.}
\end{equation}
Therefore,
the convergence~(\ref{cas6}) for\, $z\in\partial U_i$ follows immediately from~(\ref{cas5}). Take arbitrary $\e>0$ and consider the sufficiently large  radius $R_\e>R_*$ such that
\begin{equation}
\label{cas6-p}
\bigl|\gamma(z)-\frac12\bigr|<\e/2\qquad\mbox{ if \ $z\in \partial U_i$ and \ $|z|\ge R_\e$}.
\end{equation}
Since $\omega(z)>0$ in $U_i$ and $\omega(z)=0$ on $S_{R_\e}\cap \partial U_i$,
we deduce from (\ref{cas6-p}),   by continuity of $\gamma$ and by compactness arguments, that there exists $\delta=\delta_\e>0$ satisfying the condition
\begin{equation}\label{cas6-pp}|\gamma(z)-\frac12|<\e/2\qquad\mbox{ if \ $z\in  U_i$, \ $|z|=R_\e$\,, \ and \ $\omega(z)<\delta$}.
\end{equation}
Now take $R_2>R_\e$ such that
\begin{equation}\label{cas6-ppp}\omega(z)<\delta\qquad\mbox{ if $z\in U_i$ and $|z|>R_2$}.
\end{equation}
Consider an~arbitrary point $z_0\in U_i$ \, with \,$|z_0|>R_2$. Since $\omega$ is an~analytical nonconstant function,
by the classical Morse--Sard theorem on critical values  and by continuity of $\gamma$,
there exists $z_1\in U_i$ such that
\begin{equation}\label{cas6-ppp1}|z_1|>R_2, \qquad |\gamma(z_1)-\gamma(z_0)|<\frac\e2
\end{equation}
and
\begin{equation}\label{cas6-p1}\nabla\omega(z)\ne0\qquad\mbox{ if  \ $\omega(z)=\omega(z_1)$ \ and $z\in U_i$}.
\end{equation}
Denote $t_1=\omega(z_1)$,  then the connected component $L$  of the level set $\{z\in U_i:\omega(z)=t_1\}$ containing the~point $z_1$, is a smooth curve homeomorphic to the open interval $(-1,1)$ \ (indeed, this curve could not be closed because of  maximum principle for the vorticity~$\omega$\,). Evidently,
the intersection of the~curve $L$ with the circle $S_{R_\e}=\{z:|z|=R_\e\}$ \,contains at least  two points $A$ and $B$ such that $z_1$ lies between $A$ and $B$ with respect to~$L$.\footnote{Indeed, take an~arbitrary diffeomorphic parametrization
$f:(-1,1)\to L$. Then $f(s_1)=z_1$ for some $s_1\in (-1,1)$,  further, by construction we have
 \begin{equation} \label{cas6-p2}\omega(f(s))\equiv t_1>0.
\end{equation} Then the closure of~$L$ is a compact set and, of course,
\begin{equation}\label{exp-ap01}\mbox{$\dist(f(s),\partial U_i)\to 0$ as $|s|\to1$}.
\end{equation}
The~property~(\ref{cas6-p2}) guaranties  that $L$ is separated from the closed set $\{z\in\partial U_i:|z|\ge R_\e\}$. \,Therefore, by~(\ref{exp-ap01}) we have $|f(s)|<R_\e$ when $|s|$ is sufficiently close to~$1$, and this (together with the assumption~$|f(s_1)|=|z_1|>R_\e$\,) implies the existence of $s',s''\in (-1,1)$ such that
$s'<s_1<s''$ and $|f(s')|=|f(s'')|=R_\e$. Now we can take $A=f(s')$ and $B=f(s'')$.} \,By construction, $\omega|_L\equiv t_1<\delta$, \,thus by (\ref{cas6-pp}) we have
\begin{equation}\label{cas6-pp7}|\gamma(A)-\frac12|<\e/2,\qquad |\gamma(B)-\frac12|<\e/2.
\end{equation}
This implies, by virtue of the monotonicity of $\gamma$ along the curve ~$L$ \,(see~(\ref{exp-app2})\,),\, that $|\gamma(z_1)-\frac12|<\e/2$. Taking into account the second  inequality in~(\ref{cas6-ppp1}), we obtain
\begin{equation}\label{cas6-pp12}\bigl|\gamma(z_0)-\frac12\bigr|<\e.
\end{equation}
In other words, for every  point $z_0\in U_i$ with $|z_0|>R_2$ we proved  the estimate~(\ref{cas6-pp12}). Since $\e>0$ was arbitrary, the required
convergence~(\ref{cas6}) is established. \hfill $\qed$

\medskip
Since there exist only finitely many components~$U_i$, from Theorem~\ref{ttt3} we obtain immediately

\begin{cor}\label{ttt4}{\sl  The  convergence
\begin{equation}\label{cas7}\gamma(z)\to \frac12|\u_\infty|^2\qquad\mbox{ uniformly as \ }|z|\to\infty
\end{equation}
holds.
}
\end{cor}

The function $\gamma=\Phi-\omega\psi $ is closely related to $\Phi$; in particular, $\gamma=\Phi$ \ if \ $\omega=0$ or $\psi=0$. Having this in mind, it is possible
to prove the same convergence as~(\ref{cas7}) for $\Phi$ instead of $\gamma$.

We  assume without loss of generality that
\begin{equation}\label{cas7i}\u_\infty=(1,0).
\end{equation}

Recall that by
 D.~Gilbarg \& H. Weinberger results \cite{GW2} the convergence
 \begin{equation}
\label{GGWW-as1}
\displaystyle\lim_{r\to+\infty}\into0^{2\pi}|\u(r,\theta)-\u_\infty|^2d\theta=0\end{equation}
holds. In other words, since $\nabla\psi=\u^\bot=(-v,u)$, we have
 \begin{equation}
\label{GGWW-as2}
\displaystyle\lim_{r\to+\infty}\frac1r\int\limits_{|z|=r}|\nabla\psi(z)-(0,1)|^2\,ds=0.\end{equation}
 Form this fact and from the finiteness of the Dirichlet integral $\int\limits_{\Omega}|\nabla\ue|^2<\infty$ we obtain (see \cite[pages~99--100]{Amick} for details)
the following asymptotic behaviour of the stream function $\psi$:\footnote{Stream function $\psi$ is well defined by identity $\nabla\psi=\ue^\bot$ in every simply--connected subdomain of~$\Omega_*$; in particular,
$\psi$ is well-defined in  intersection of~$\Omega_*$ with every of the four half spaces $\{(x,y)\in\R^2: x>0\}$, \ $\{(x,y)\in\R^2: x\le 0\}$, \ $\{(x,y)\in\R^2: y>0\}$, \ $\{(x,y)\in\R^2: y\le0\}$. Since these definitions of~$\psi$ differ only by some additive constants, they have
no influence on the asymptotic properties discussed  here.}
 \begin{equation}
\label{GGWW-as3}
\displaystyle\lim_{r\to+\infty}\frac1r|\psi(x,y)-y|=0,\end{equation}
where $r=\sqrt{x^2+y^2}$. For any $\alpha>0$ denote by $\SA$ the sector
$$
\SA=\{z=(x,y)\in\Omega_*:\frac{|y|}{|x|}\ge\alpha\}.
$$
Since
$r\le c_\alpha |y|$ for $z\in \SA$, from $(\ref{GGWW-as3})$  it follows that
 \begin{equation}
\label{GGWW-as4}
\displaystyle\lim_{(x,y)\in S_\alpha, \ \sqrt{x^2+y^2}\to\infty}\biggl|\frac{\psi(x,y)}y-1\biggr|=0.
\end{equation}
Let us prove the convergence of~$\Phi$ in any sector $\SA$.

\begin{lem}[see Theorem~17 and Corollary 18 on  page 101 in~\cite{Amick}]\label{ttt5}{\sl  For any $\alpha>0$ the uniform convergences
\begin{equation}\label{cas11}
|z|\omega(z)\to0 \qquad\mbox{ as }|z|\to\infty, \ z\in \SA,
\end{equation}
\begin{equation}\label{cas12}\Phi(z)\to  \frac12|\u_\infty|^2\qquad\mbox{ as }|z|\to\infty, \ z\in \SA.
\end{equation}
hold.
}
\end{lem}

{\sc Proof}.
Fix $\alpha>0$. Then
\begin{equation}\label{cas18}\forall z=(x,y)\in\SAA:|z|\le \tilde c_\alpha |y|.
\end{equation}
Take $z_0=(x_0,y_0)\in\SA$.  Without loss of generality assume that $y_0>0$.
Since $$\int\limits_{\Omega_*}|\nabla\Phi|^2<\infty,$$ from Lemma~\ref{lt2}, from the~uniform convergence of the pressure to zero (see~(\ref{cas4})) and from average convergence of the velocity to~$\u_\infty=(1,0)$ \  (see~(\ref{GGWW-as1})\,),  we  have that
\begin{equation}\label{cas13}\exists r\in[\frac14y_0,\frac12y_0]:\qquad\sup\limits_{|z-z_0|=r}\bigl|\Phi(z)-\frac12\bigr|\le \e_1(r_0),
\end{equation}
where $r_0=|z_0|$ and $\e_1(r_0)\to0$ uniformly as $r_0\to\infty$ (of course, this function $\e_1(r_0)$ depends  also on the parameter~$\alpha$ fixed above).

From  \eqref{cas13} and from Corollary~\ref{ttt4} we have
\begin{equation}\label{cas13o}\sup\limits_{|z-z_0|=r}|\omega(z)\psi(z)|\le \e_2(r_0),
\end{equation}
where again $\e_2(r_0)\to0$ uniformly as $r_0\to\infty$.
Denote by $B_0$ the disk
$\{z\in \R^2:|z-z_0|\le r\}$.  By construction,
$$B_0\subset \SAA.$$
Then by (\ref{GGWW-as4}),
\begin{equation}
\label{GGWW-as9}
\sup_{(x,y)\in B_0}\biggl|\frac{\psi(x,y)}y-1\biggr|\to0\mbox{ as } r_0\to\infty.\end{equation}
In particular,
\begin{equation}
\label{GGWW-as10}
\psi(y)\ge d_\alpha r_0\end{equation}
if $r_0$ is sufficiently large, here the constant~$d_\alpha$ depends  on~$\alpha$ only. From (\ref{GGWW-as10}) and  (\ref{cas13o}) we obtain immediately that
\begin{equation}\label{cas21}\sup\limits_{|z-z_0|=r}|\omega(z)|\le \frac1{r_0}\e_3(r_0),
\end{equation}
where again $\e_3(r_0)\to0$ uniformly as $r_0\to\infty$.
By maximum principle,
\begin{equation}\label{cas22}|\omega(z_0)|\le \frac1{r_0}\e_3(r_0).
\end{equation}
Thus, we have proved the asymptotic estimate~(\ref{cas11}). Then the convergence (\ref{cas12}) follows immediately from~(\ref{cas11}) and  (\ref{cas7}). 
\hfill$\qed$

\medskip

The convergence of~$\Phi$ {\it outside} of the sectors $\SA$ is more delicate and subtle question. Ch. Amick solved this problem~\cite{Amick} using level sets of the stream function~$\psi$.

Define the stream function in the half-domain~$\Omega_+=\{(x,y):x\ge0, \ x^2+y^2\ge R_*^2\}$ and consider the set~$C_+=\{z\in\Omega_+:\psi(z)=0\}$
\footnote{The asymptotic behavior of~$\psi(x,y)$ is similar to that of the linear function~$g(x,y)=y$.  Since the level set $\{(x,y)\in\Omega_+:g(x,y)=0\}$ is a ray $\{(x,y)\in\Omega_+: 
y=0\}$, the set $C_+$  goes to infinity as well, see also Lemma~\ref{ttt7} for the precise formulation.}.
Then $\gamma=\Phi$ on \,$C_+$ and from the convergence of~$\gamma$ (\ref{cas7}) we obtain immediately that $\frac12|\nabla\psi(z)|^2=\frac12|\u(z)|^2\to\frac12$ when $|z|\to\infty$, $z\in C_+$. In particular,
$\nabla\psi\ne0$ on $C_+$ if we choose the parameter $R_*$ sufficiently large. Using similar arguments, Amick proved that
  the set~$C_+$ has very simple geometrical structure.

\begin{lem}[see Lemma 20 on  page 104 in~\cite{Amick}]\label{ttt7}{\sl  If the number \ $R_*$ is chosen large enough, then the set $C_+$ is a smooth curve
$$C_+=\bigl\{(p_+(s),q_+(s)):s\in[0,+\infty)\bigr\},$$
here $p_+$ and $q_+$ are real-analytic functions on~$[0,\infty)$, \ $p_+(s)\to\infty$ and $\frac{q_+(s)}{p_+(s)}\to0$ as $s\to\infty$. In addition,
\begin{equation}\label{cas23}|\u(p_+(s),q_+(s))|\to|\u_\infty|\quad\mbox{ as }s\to\infty.
\end{equation}
}
\end{lem}

Of course, the similar assertion holds for another  half-domain~$\Omega_-=\{(x,y):x\le0, \ x^2+y^2\ge R_*^2\}$.  Using this Lemma and some classical estimates for the Laplace operator
(recall, that $\omega=\Delta\psi$\,), Amick proved the required assertion:

\begin{ttt}[see Theorem 21~(a) on  page 1045 in~\cite{Amick}]\label{ttt8}{\sl  The convergence \begin{equation}
\label{2c-asp---}
|\u(z)|\to|\u_\infty|\qquad\mbox{ uniformly as \ }|z|\to\infty.
\end{equation}holds. }
\end{ttt}

\begin{rem}\label{rem20}{\rm The proof of Theorem~\ref{ttt8} could be essentially simplified in comparison with the original version of~\cite{Amick}. Indeed, from the convergence~(\ref{cas23}) on the curve~$C_+$, using the Lemmas~\ref{lt1}--\ref{lt2} it is very easy to derive that there exists $\sigma>0$ such that for any $z\in C_+$ with sufficiently large value~$|z|$ we have
$$\biggl|\frac1{r}\int\limits_{|\xi-z|=r}\u(\xi)\,ds\biggr|>\sigma$$
for all $r\in(0, \frac45|z|]$. Then the arguments of the proof of  Lemma~\ref{T3}~(ii) of the present paper give us that
\begin{equation}\label{cas24}\u(p_+(s),q_+(s))\to\u_\infty\quad\mbox{ as }s\to\infty
\end{equation}
instead of~(\ref{cas23}). This more strong convergence allows to simplify some technical moments in the proof of~\cite[Theorem 21~(a)]{Amick}, see also \cite[Theorem 21~(c)]{Amick}.}
\end{rem}

\

 {\em {\bf Acknowledgment.}} M.~Korobkov was partially supported by
the Ministry of Education and Science of the Russian Federation (the Project number 1.3087.2017/4.6) and by the Russian Federation for Basic Research  (Project no.18-01-00649a).

The research of K. Pileckas was funded by the grant No. S-MIP-17-68 from the Research Council of Lithuania.

\medskip 

 { {\bf  Conflict of interest:}} The authors declare that they have no conflict of interest.

{\small 
\end{document}